\input ./style/arxiv-vmsta.cfg
\documentclass[numbers,compress]{vmsta}
\usepackage{enumitem}
\usepackage{amssymb}
\usepackage{mathscinet}
\usepackage{vtexurl}

\volume{2}
\pubyear{2015}
\firstpage{1}
\lastpage{15}
\doi{10.15559/15-VMSTA20}

\startlocaldefs

\newtheorem{thm}{Theorem}

\newtheorem{lemma}{Lemma}

\theoremstyle{definition}
\newtheorem{remark}{Remark}

\newtheorem{assump}{Assumption}

\def\ex{{\rm{\mathbb E\,}}}

\newcommand{\pp}{\mathbb{P}}
\newcommand{\ee}{\mathbb{E}}

\newcommand{\bb}[1]{\mathbb{ #1}}
\newcommand{\RR}{\mathbb{R}}
\newcommand{\dd}{\mathrm{d}}
\def\dd{\mathrm{d}}
\newcommand{\eps}{\varepsilon}

\allowdisplaybreaks

\newcommand{\rrvert}{\vert}
\newcommand{\llvert}{\vert}
\allowdisplaybreaks
\endlocaldefs

\begin{document}
\begin{frontmatter}

\title{Nonparametric Bayesian inference for multidimensional compound
Poisson processes}

\author[a]{\inits{S.}\fnm{Shota}\snm{Gugushvili}\corref{cor1}}\email
{shota.gugushvili@math.leidenuniv.nl}
\cortext[cor1]{Corresponding author.}

\author[b]{\inits{F.}\fnm{Frank}\snm{van der Meulen}}\email
{f.h.vandermeulen@tudelft.nl}

\author[c]{\inits{P.}\fnm{Peter}\snm{Spreij}}\email{spreij@uva.nl}

\address[a]{Mathematical Institute, Leiden University, P.O. Box 9512,\\
2300 RA Leiden, The~Netherlands}
\address[b]{Delft Institute of Applied Mathematics, Faculty of
Electrical Engineering, Mathematics and Computer Science, Delft
University of Technology, Mekelweg 4, 2628 CD Delft, The Netherlands}
\address[c]{Korteweg--de Vries Institute for Mathematics, University of
Amsterdam,\\ P.O.\ Box 94248, 1090 GE Amsterdam, The Netherlands}

\markboth{S. Gugushvili et al.}{Nonparametric Bayesian inference for multidimensional compound
Poisson processes}

\begin{abstract}
Given a sample from a discretely observed multidimensional compound
Poisson process, we study the problem of nonparametric estimation of
its jump size density $r_0$ and intensity $\lambda_0$. We take a
nonparametric Bayesian approach to the problem and determine posterior
contraction rates in this context, which, under some assumptions, we
argue to be optimal posterior contraction rates. In particular, our
results imply the existence of Bayesian point estimates that converge
to the true parameter pair $(r_0,\lambda_0)$ at these rates. To the
best of our knowledge, construction of nonparametric density estimators
for inference in the class of discretely observed multidimensional
L\'evy processes, and the study of their rates of convergence is a~new
contribution to the literature.
\end{abstract}

\received{24 December 2014}
%
\revised{27 February 2015}
%
\accepted{1 March 2015}
\publishedonline{13 March 2015}
\end{frontmatter}

\section{Introduction}
\label{intro}

Let $N=(N_t)_{t\geq0}$ be a Poisson process of constant intensity
$\lambda>0$, and let $\{Y_j\}$ be independent and identically
distributed (i.i.d.) $\RR^d$-valued random vectors defined on the same
probability space and having a common distribution function $R$, which
is assumed to be absolutely continuous with respect to the Lebesgue
measure with density $r$. Assume that $N$ and $\{Y_j\}$ are independent
and define the $\mathbb{R}^d$-valued process $X=(X_t)_{t\geq0}$ by
\begin{equation*}
X_t=\sum_{j=1}^{N_t}Y_j.
\end{equation*}
The process $X$ is called a compound Poisson process (CPP) and forms a
basic stochastic model in a variety of applied fields, such as, for
example,\ risk theory and queueing; see \cite{embrechts97,prabhu98}.

Suppose that, corresponding to the true parameter pair $(\lambda
_0,r_0)$, a sample $X_{\varDelta}$, $X_{2\varDelta},\ldots,X_{n\varDelta}$ from
$X$ is available, where the sampling mesh $\varDelta>0$ is assumed to be
fixed and thus independent of $n$. The problem we study in this note is
nonparametric estimation of $r_0$ (and of $\lambda_0$). This is
referred to as decompounding and is well studied for one-dimensional
CPPs; see \cite{buchmann03,buchmann04,comte13,duval12,vanes07}.
Some practical situations in which this problem may
arise are listed in \cite[p.\ 3964]{duval12}. However, the methods used
in the above papers do not seem to admit (with the exception of
\cite{vanes07}) a generalization to the multidimensional setup. This is
also true for papers studying nonparametric inference for more general
classes of L\'evy processes (of which CPPs form a~particular class),
such as, for example,\ \cite{comte10,comte11,neumann09}. In
fact, there is a dearth of publications dealing with nonparametric
inference for multidimensional L\'evy processes. An exception is
\cite{bucher13}, where the setup is however specific in that it is
geared to inference in L\'evy copula models and that, unlike the
present work, the high-frequency sampling scheme is assumed
($\varDelta=\varDelta_n\rightarrow0$ and $n\varDelta _n\rightarrow\infty$).

In this work, we will establish the posterior contraction rate in a
suitable metric around the true parameter pair $(\lambda_0,r_0)$. This
concerns study of asymptotic frequentist properties of Bayesian
procedures, which has lately received considerable attention in the
literature (see, e.g.,\  \cite{ghosal01,ghosal00}), and is
useful in that it provides their justification from the frequentist
point of view. Our main result says that for a $\beta$-H\"older regular
density $r_0$, under some suitable additional assumptions on the model
and the prior, the posterior contracts at the rate $n^{-\beta/(2\beta
+d)}(\log n)^{\ell}$, which, perhaps up to a logarithmic factor, is
arguably the optimal posterior contraction rate in our problem.
Finally, our Bayesian procedure is adaptive: the construction of our
prior does not require knowledge of the smoothness level $\beta$ in
order to achieve the posterior contraction rate given above.

The proof of our main theorem employs certain results from
\cite{ghosal01,shen13} but involves a substantial number of
technicalities specifically characteristic of decompounding.

We remark that a practical implementation of the Bayesian approach to
decompounding lies outside the scope of the present paper. Preliminary
investigations and a small scale simulation study we performed show
that it is feasible and under certain conditions leads to good results.
However, the technical complications one has to deal with are quite
formidable, and therefore the results of our study of implementational
aspects of decompounding will be reported elsewhere.

The rest of the paper is organized as follows. In the next section, we
introduce some notation and recall a number of notions useful for our
purposes. Section \ref{main} contains our main result, Theorem \ref
{mainthm}, and a brief discussion on it. The proof of Theorem \ref
{mainthm} is given in Section \ref{proofs}. Finally, Section \ref
{pr.lem.1} contains the proof of the key technical lemma used in our proofs.

\section{Preliminaries}

Assume without loss of generality that $\varDelta=1$, and let
$Z_{i}=X_{i}-X_{i-1}$, $i=1,\ldots,n$. The $\mathbb{R}^d$-valued random
vectors ${Z_i}$ are i.i.d.\ copies of a random vector
\begin{equation*}
Z=\sum_{j=1}^T Y_j,
\end{equation*}
where $\{Y_j\}$ are i.i.d.\ with distribution function $R_0$, whereas
$T$, which is independent of $\{Y_j\}$, has the Poisson distribution
with parameter $\lambda_0$. The problem of decompounding the jump size
density $r_0$ introduced in Section \ref{intro} is equivalent to
estimation of $r_0$ from observations ${\mathcal{Z}_n}=\{Z_1,Z_2,\ldots
,Z_n\}$, and we will henceforth concentrate on this alternative
formulation. We will use the following notation:
\begin{description}
\item[$\mathbb{P}_r$] law of \xch{$Y_1$,}{$Y_1$}
\item[$\mathbb{Q}_{\lambda, r}$] law of \xch{$Z_1$,}{$Z_1$}
\item[$\mathbb{R}_{\lambda, r}$] law of \xch{$X=(X_t,\, t\in[0,1])$.}{$X=(X_t,\, t\in[0,1])$}
\end{description}

\subsection{Likelihood}
We will first specify the dominating measure for $\mathbb{Q}_{\lambda,
r}$, which allows us to write down the likelihood in our model. Define
the random measure $\mu$ by
\begin{equation*}
\mu( B ) = \bigl\{ \# t :(t,X_t-X_{t-})\in B \bigr\},
\quad B \in\mathcal {B}\bigl([0,1]\bigr)\otimes\mathcal{B}\bigl(
\mathbb{R}^d\setminus\{0\}\bigr).
\end{equation*}
Under $\mathbb{R}_{\lambda,r}$, the random measure $\mu$ is a Poisson
point process on $[0,1]\times(\mathbb{R}^d\setminus\{0\})$ with
intensity measure $ {\varLambda}(\dd t,dx)={\lambda} \dd t {r}(x) \dd x$.
Provided that $\lambda,\widetilde{\lambda}>0$, and $\widetilde{r}>0$,
by formula (46.1) on p.~262 in~\cite{skorohod64} we have
%
\begin{equation}
\label{eq:contlik} \frac{\mathrm{d} \mathbb{R}_{\lambda,r}}{ \mathrm{d}\mathbb
{R}_{\widetilde{\lambda},\widetilde{r}} }(X)=\exp \Biggl( \int_0^1
\int_{\mathbb{R}^d} \log \biggl( \frac{\lambda r(x)}{{\widetilde{\lambda}}
\widetilde{r}(x)} \biggr) \mu(\dd t,
\dd x) - (\lambda-{\widetilde {\lambda}}) \Biggr).
\end{equation}
The density $k_{\lambda,r}$ of $\mathbb{Q}_{\lambda,r}$ with respect to
$\mathbb{Q}_{\widetilde{\lambda},\widetilde{r}}$ is then given by the
conditional expectation
%
\begin{equation}
\label{eq:k} k_{\lambda,r}(x)=\ex_{\widetilde{\lambda},\widetilde{r}} \biggl(
\frac
{\mathrm{d} \mathbb{R}_{\lambda,r}}{ \mathrm{d}\mathbb{R}_{\widetilde
{\lambda},\widetilde{r}} }(X)\, \Big|\, X_1=x \biggr),
\end{equation}
where the subscript in the conditional expectation operator signifies
the fact that it is evaluated under $\mathbb{R}_{ \widetilde{\lambda
},\widetilde{r} }$; see Theorem 2 on p.~245 in~\cite{skorohod64} and
Corollary 2 on p.~246 there. Hence, the likelihood (in the parameter
pair $(\lambda,r)$) associated with the sample $\mathcal{Z}_n$ is given
by
%
\begin{equation}
\label{ln} L_n(\lambda,r)=\prod_{i=1}^n
k_{\lambda,r}(Z_i).
\end{equation}

\subsection{Prior}
\label{prior.sec}
We will use the product prior $\varPi=\varPi_1\times\varPi_2$ for $(\lambda
_0,r_0)$. The prior $\varPi_1$ for $\lambda_0$ will be assumed to be
supported on the interval $[\underline{\lambda},\overline{\lambda}]$
and to possess a density $\pi_1$ with respect to the Lebesgue measure.

The prior for $r_0$ will be specified as a Dirichlet process mixture of
normal densities. Namely, introduce a convolution density
%
\begin{equation}
\label{fjs} r_{F,\varSigma}(x)=\int\phi_{\varSigma}(x-z)F(\dd z),
\end{equation}
where $F$ is a distribution function on $\mathbb{R}^d$, $\varSigma$ is a
$d\times d$ positive definite real matrix, and $\phi_{\varSigma}$ denotes
the density of the centered $d$-dimensional normal distribution with
covariance matrix $\varSigma$. Let $\alpha$ be a finite measure on $\RR
^d$, and let $\mathcal{D}_{\alpha}$ denote the Dirichlet process
distribution with base measure $\alpha$ (see \cite{ferguson73} or,
alternatively, \cite{ghosal10} for a modern overview). Recall that if
$F\sim\mathcal{D}_{\alpha}$, then for any Borel-measurable partition
$B_1,\ldots,B_k$ of $\RR^d$, the distribution of the vector
$(F(B_1),\ldots,F(B_k))$ is the $k$-dimensional Dirichlet distribution
with parameters $\alpha(B_1),\ldots,\alpha(B_k)$. The Dirichlet process
location mixture of normals prior $\varPi_2$ is obtained as the law of the
random function $r_{F,\varSigma}$, where $F\sim\mathcal{D}_{\alpha}$ and
$\varSigma\sim G$ for some prior distribution function $G$ on the set of
$d\times d$ positive definite matrices. For additional information on
Dirichlet process mixtures of normal densities, see, for example,\ the
original papers \cite{ferguson83} and \cite{lo84}, or a recent paper
\cite{shen13} and the references therein.

\subsection{Posterior}
Let $\mathcal{R}$ denote the class of probability densities of the form
\eqref{fjs}. By Bayes' theorem, the posterior measure of any measurable set
$A\subset(0,\infty)\times\mathcal{R}$
is given by
\begin{equation*}
\varPi(A|\mathcal{Z}_n)=\frac{ \iint_A L_n(\lambda,r) \dd\varPi_1(\lambda)
\mathrm{d}\varPi_2(r) }{ \iint L_n(\lambda,r) \dd\varPi_1(\lambda) \mathrm
{d}\varPi_2(r) }.
\end{equation*}
The priors $\varPi_1$ and $\varPi_2$ indirectly induce the prior $\varPi= \varPi
_1\times\varPi_2$ on the collection of densities $k_{\lambda,r}$. We will
use the symbol $\varPi$ to signify both the prior on $(\lambda_0,r_0)$ and
the density $k_{\lambda_0,r_0}$. The posterior in the first case will
be understood as the posterior for the pair $(\lambda_0,r_0)$, whereas
in the second case as the posterior for the density $k_{\lambda_0,r_0}$.
Thus, setting $\overline{A}=\{ k_{\lambda,r}:(\lambda,r)\in A \}$, we have
\begin{equation*}
\varPi(\overline{A}|\mathcal{Z}_n)=\frac{ \int_{\overline A} L_n(k) \mathrm
{d}\varPi(k) }{ \int L_n(k) \mathrm{d}\varPi(k) }.
\end{equation*}
In the Bayesian paradigm, the posterior encapsulates all the
inferential conclusions for the problem at hand. Once the posterior is
available, one can next proceed with computation of other quantities of
interest in Bayesian statistics, such as Bayes point estimates or
credible sets.

\begingroup\abovedisplayskip=8pt\belowdisplayskip=8pt
\subsection{Distances}
The Hellinger distance $h(\mathbb{Q}_{0},\mathbb{Q}_{1})$ between two
probability laws $\mathbb{Q}_{0}$ and $\mathbb{Q}_{1}$ on a measurable
space $(\varOmega,\mathfrak{F})$ is given by
\[
h(\mathbb{Q}_{0},\mathbb{Q}_{1})= \biggl( \int \bigl(
\mathrm{d}\mathbb {Q}_{0}^{1/2} - \mathrm{d}
\mathbb{Q}_{1}^{1/2} \bigr)^2
\biggr)^{1/2}.
\]
Assuming that $\mathbb{Q}_0\ll\mathbb{Q}_1$, the Kullback--Leibler
divergence $\mathrm{K}(\mathbb{Q}_{0},\mathbb{Q}_{1})$ is
\[
\mathrm{K}(\mathbb{Q}_{0},\mathbb{Q}_{1})= \int\log
\biggl(\frac{
\mathrm{d}\mathbb{Q}_{0}}{ \mathrm{d}\mathbb{Q}_{1} } \biggr) \mathrm {d}\mathbb{Q}_{0}.
\]
We also define the $\mathrm{V}$-discrepancy by
\[
\mathrm{V}(\mathbb{Q}_{0},\mathbb{Q}_{1}) =\int
\log^2 \biggl(\frac{ \mathrm{d}\mathbb{Q}_{0}}{ \mathrm{d}\mathbb
{Q}_{1} } \biggr) \mathrm{d}
\mathbb{Q}_{0}.
\]
In addition, for positive real numbers $x$ and $y$, we put
\begin{align*}
\mathrm{K}(x,y) & = x\log\frac{x}{y}-x+y,
\\
\mathrm{V}(x,y) & = x\log^2\frac{x}{y},
\\
h(x,y) & = \bigl|\sqrt{x}-\sqrt{y}\bigr|.
\end{align*}
Using the same symbols $\mathrm{K}$, $\mathrm{V}$, and $h$ is justified
as follows. Suppose that $\varOmega$ is a singleton $\{\omega\}$ and
consider the Dirac measures $\delta_x$ and $\delta_y$ that put masses
$x$ and $y$, respectively, on $\varOmega$. Then $\mathrm{K}(\delta_x,\delta
_y)=\mathrm{K}(x,y)$, and similar equalities are valid for the $\mathrm
{V}$-discrepancy and the Hellinger distance.

\subsection{Class of locally $\beta$-H\"older functions}

For any $\beta\in\RR$, by $\lfloor\beta\rfloor$ we denote the largest
integer strictly smaller than $\beta$, by $\mathbb{N}$ the set of
natural numbers, whereas $\mathbb{N}_0$ stands for the union $\mathbb
{N}\xch{\cup}{\bigcup}\{0\}$. For a multiindex $k=(k_1,\ldots,k_d)\in\mathbb
{N}_0^d$, we set $k_{.}=\sum_{i=1}^d k_i$. The usual Euclidean norm of
a vector $y\in\RR^d$ is denoted by $\|y\|$.

Let $\beta>0$ and $\tau_0\geq0$ be constants, and let $L:\RR
^d\rightarrow\RR_{+}$ be a measurable function. We define the class
$\mathcal{C}^{\beta,L,\tau_0}(\RR^d)$ of locally $\beta$-H\"older
regular functions as the set of all functions $r:\RR^d\rightarrow\RR$
such that all mixed partial derivatives $D^k r$ of $r$ up to order
$k_{.}\leq\lfloor\beta\rfloor$ exist and, for every $k$ with
$k_{.}=\lfloor\beta\rfloor$, satisfy
\begin{equation*}
\bigl| \bigl(D^k r\bigr) (x+y) - \bigl(D^k\bigr)r(x) \bigr|\leq
L(x)\exp\bigl(\tau_0 \|y\|^2\bigr)\|y\|^{\beta
-\lfloor\beta\rfloor},
\quad x,y\in\RR^d.
\end{equation*}
See p.~625 in \cite{shen13} for this class of functions.

\section{Main result}
\label{main}

Define the complements of the Hellinger-type neighborhoods of $(\lambda
_0,r_0)$ by
\begin{equation*}
A(\varepsilon_n,M)=\bigl\{ (\lambda,r) : h( \mathbb{Q}_{\lambda
_0,r_0},
\mathbb{Q}_{\lambda,r} )> M \varepsilon_n \bigr\},
\end{equation*}
where $\{\varepsilon_n\}$ is a sequence of positive numbers. We say that
$\varepsilon_n$ is a posterior contraction rate if there exists a
constant $M>0$ such that
\begin{equation*}
\varPi\bigl(A(\varepsilon_n,M)\big|\mathcal{Z}_n\bigr)
\rightarrow0
\end{equation*}
as $n\rightarrow\infty$ in $\mathbb{Q}_{\lambda_0,r_0}^n$-probability.

The $\varepsilon$-covering number of a subset $B$ of a metric space
equipped with the metric $\rho$ is the minimum number of $\rho$-balls
of radius $\varepsilon$ needed to cover it. Let
$\mathcal{Q}$ be a set of CPP laws $\bb{Q}_{\lambda,r}$. Furthermore,
we set
%
\begin{equation}
\label{set_b} B(\varepsilon,\mathbb{Q}_{\lambda_0,r_0})=\bigl\{ (\lambda,r):
\mathrm {K}(\mathbb{Q}_{\lambda_0,r_0},\mathbb{Q}_{\lambda,r})\leq\varepsilon
^2, \mathrm{V}(\mathbb{Q}_{\lambda_0,r_0},\mathbb{Q}_{\lambda,r})
\leq \varepsilon^2\bigr\}.
\end{equation}
We recall the following general result on posterior contraction rates.

\begin{thm}[\cite{ghosal01}]
\label{thm2.1ghosal01} Suppose that for positive sequences
$\overline{\varepsilon}_n,\widetilde {\varepsilon}_n\rightarrow0$ such
that\linebreak $n\min(\overline{\varepsilon
}_n^2,\widetilde{\varepsilon}_n^2)\rightarrow\infty$, constants
$c_1,c_2,c_3,c_4>0$, and sets $\mathcal{Q}_n\subset\mathcal{Q}$, we
have
\begin{align}
\log N(\overline{\varepsilon}_n,\mathcal{Q}_n,h) & \leq
c_1 n\overline {\varepsilon}_n^2,\label{c1}
\\
\varPi( \mathcal{Q}\setminus\mathcal{Q}_n ) & \leq c_3
e^{- n\widetilde
{\varepsilon}_n^2 (c_2+4) }, \label{c2}
\\
\varPi\bigl( B(\widetilde{\varepsilon}_n,\mathbb{Q}_{\lambda_0,r_0})
\bigr) & \geq c_4 e ^{-c_2 n\widetilde{\varepsilon}_n^2} \label{c3}.
\end{align}
Then, for $\varepsilon_n=\max(\overline{\varepsilon}_n,\widetilde
{\varepsilon}_n)$ and a constant $M>0$ large enough, we have that
%
\begin{equation}
\label{postrate1} \varPi\bigl( A(\varepsilon_n,M)\big|
\mathcal{Z}_{n} \bigr) \rightarrow0
\end{equation}
as $n\rightarrow\infty$ in $\mathbb{Q}_{\lambda_0,r_0}^n$-probability,
assuming that the i.i.d.\ observations $\{Z_j\}$ have been generated
according to $\mathbb{Q}_{\lambda_0,r_0}$. 
\end{thm}

In order to derive the posterior contraction rate in our problem, we
impose the following conditions on the true parameter pair $(\lambda_0,r_0)$.
\begin{assump}\label{ass:truth}
Denote by $(\lambda_0, r_0)$ the true parameter values for the compound
Poisson process.
\begin{enumerate}[label={\rm(\roman*)}]
\item$\lambda_0$ is in a compact set $[\underline{\lambda}, \overline
{\lambda}]\subset(0,\infty)$;
\item
The true density $r_0$ is bounded, belongs to the set $\mathcal
{C}^{\beta,L,\tau_0}(\RR^d)$, and additionally satisfies, for some
$\varepsilon>0$ and all $k\in\mathbb{N}_0^d,\, k_{.}\leq\beta$,
\begin{equation*}
\int \biggl( \frac{L}{r_0} \biggr)^{(2\beta+\varepsilon)/\beta} r_0<\infty,
\qquad \int \biggl(\frac{|D^k r_0|}{r_0} \biggr)^{(2\beta+\varepsilon)/k} r_0<
\infty.
\end{equation*}
Furthermore, we assume that there exist strictly positive constants
$a,b,c$, and $\tau$ such that
\begin{equation*}
r_0(x) \leq c \exp\bigl(-b\|x\|^{\tau}\bigr), \quad\|x\|>a.
\end{equation*}
\end{enumerate}
\end{assump}
The conditions on $r_0$ come from Theorem 1 in \cite{shen13} and are
quite reasonable. They simplify greatly when $r_0$ has a compact support.

We also need to make some assumptions on the prior $\varPi$ defined in
Section~\ref{prior.sec}.
\begin{assump}\label{ass:prior}
The prior $\varPi=\varPi_1\times\varPi_2$ on $(\lambda_0,r_0)$ satisfies the
following assumptions:
\begin{enumerate}[label={\rm(\roman*)}]
\item The prior $\varPi_1$ on $\lambda$ has a density $\pi_1$ (with
respect to the Lebesgue measure) that is supported on the finite
interval $[\underline{\lambda},\overline{\lambda}]\subset(0,\infty)$
and is such that
%
\begin{equation}
\label{pi1} 0<\underline{\pi}_1 \leq\pi_1(\lambda)
\leq\overline{\pi}_1<\infty, \quad\lambda\in[\underline{\lambda},
\overline{\lambda}],
\end{equation}
for some constants $\underline{\pi}_1$ and $\overline{\pi}_1$;

\item The base measure $\alpha$ of the Dirichlet process prior $\mathcal
{D}_{\alpha}$ is finite and possesses a strictly positive density on
$\RR^d$ such that for all sufficiently large $x>0$ and some strictly
positive constants $a_1,b_1$, and $C_1$,
\begin{equation*}
1-\overline{\alpha}\bigl([-x,x]^d\bigr)\leq b_1 \exp
\bigl(-C_1 x^{a_1}\bigr),
\end{equation*}
where $\overline{\alpha}(\cdot)=\alpha(\cdot)/\alpha(\RR^d)$\xch{;}{.}

\item There exist strictly positive constants $\kappa,a_2$, $a_3$,
$a_4$, $a_5$, $b_2$, $b_3$, $b_4$, $C_2$, $C_3$ such that for all $x>0$
large enough,
\begin{equation*}
G\bigl( \varSigma:\operatorname{eig}_d\bigl(\varSigma^{-1}
\bigr) \geq x \bigr) \leq b_2 \exp\bigl(-C_2
x^{a_2}\bigr),
\end{equation*}
for all $x>0$ small enough,
\begin{equation*}
G \bigl(\varSigma: \operatorname{eig}_1\bigl(\varSigma^{-1}
\bigr)<x \bigr)\leq b_3 x^{a_3},
\end{equation*}
and for any $0<s_1\leq\cdots\leq s_d$ and $t\in(0,1)$,
\begin{equation*}
G\bigl( \varSigma:s_j<\operatorname{eig}_j\bigl({
\varSigma^{-1}}\bigr)<s_j(1+t),j=1,\ldots,d \bigr)\geq
b_4 s_1^{a_4}t^{a_5}\exp
\bigl(-C_3 s_d^{\kappa/2}\bigr).
\end{equation*}
Here $\operatorname{eig}_j({\varSigma^{-1}})$ denotes the $j$th smallest
eigenvalue of the matrix $\varSigma^{-1}$.
\end{enumerate}
\end{assump}

This assumption comes from \cite[p.~626]{shen13}, to which we refer for
an additional discussion. In particular, it is shown there that an
inverse Wishart distribution (a popular prior distribution for
covariance matrices) satisfies the assumptions on $G$ with $\kappa=2$.
As far as $\alpha$ is concerned, we can take it such that its rescaled
version $\overline{\alpha}$ is a nondegenerate Gaussian distribution on
$\mathbb{R}^d$.

\begin{remark}
\label{cond_pi1}
Assumption \eqref{pi1} requiring that the prior density $\pi_1$ is
bounded away from zero on the interval $[\underline{\lambda},\overline
{\lambda}]$ can be relaxed to allowing it to take the zero value at the
end points of this interval, provided that $\lambda_0$ is an interior
point of $[\underline{\lambda},\overline{\lambda}]$.
\end{remark}

We now state our main result.

\begin{thm}
\label{mainthm}
Let Assumptions~\ref{ass:truth} and~\ref{ass:prior} hold. Then
there exists a constant $M>0$ such that, as $n\rightarrow\infty$,
\begin{equation*}
\varPi \bigl( A \bigl({(\log n)^{\ell} }
n^{-\gamma},M \bigr)
\big| \mathcal{Z}_n \bigr)\rightarrow0
\end{equation*}
in $\mathbb{Q}_{\lambda_0,r_0}^n$-probability. Here
\begin{equation*}
\gamma=\frac{\beta}{2\beta+d^*}, \qquad\ell>\ell_0=\frac{d^*(1+1/\tau
+1/\beta)+1}{2+d^*/\beta},
\qquad d^*=\max(d,\kappa).
\end{equation*}
\end{thm}

We conclude this section with a brief discussion on the obtained
result: the logarithmic factor $(\log n)^{\ell}$ is negligible for
practical purposes. If $\kappa=1$, then the posterior contraction rate
obtained in Theorem \ref{mainthm} is essentially $n^{-2\beta/(2\beta
+d)}$, which is the minimax estimation rate in a number of
nonparametric settings. This is arguably also the minimax estimation
rate in our problem as well (cf.\ Theorem 2.1 in \cite{gugu08} for a
related result in the one-dimensional setting), although here we do not
give a formal argument. Equally important is the fact that our result
is adaptive: the posterior contraction rate in Theorem \ref{mainthm} is
attained without the knowledge of the smoothness level $\beta$ being
incorporated in the construction of our prior $\varPi$. Finally, Theorem
\ref{mainthm}, in combination with Theorem 2.5 and the arguments on
pp.\ 506--507 in \cite{ghosal00}, implies the existence of Bayesian
point estimates achieving (in the frequentist sense) this convergence rate.

\begin{remark}
After completion of this work, we learned about the paper \cite{donnet14} that deals with nonparametric Bayesian estimation of
intensity functions for Aalen counting processes. Although CPPs are in
some sense similar to the latter class of processes, they are not
counting processes. An essential difference between our work and
\cite{donnet14} lies in the fact that, unlike \cite{donnet14},
ours deals with discretely observed multidimensional processes. Also
\cite{donnet14} uses the log-spline prior, or the Dirichlet
mixture of uniform densities, and not the Dirichlet mixture of normal
densities as the prior.
\end{remark}
\endgroup

\section{Proof of Theorem \ref{mainthm}}
\label{proofs}
\begingroup\abovedisplayskip=7pt\belowdisplayskip=7pt
The proof of Theorem \ref{mainthm} consists in verification of the
conditions in Theorem~\ref{thm2.1ghosal01}. The following lemma plays
the key role.

\begin{lemma}\label{lem:ineq}
The following estimates are valid:
\begin{align}
\mathrm{K}(\mathbb{Q}_{\lambda_0,r_0},\mathbb{Q}_{\lambda,r}) & \leq
\lambda_0 \mathrm{K}(\mathbb{P}_{r_0},\mathbb{P}_{r})+
\mathrm{K}(\lambda _0,\lambda),\label{eq:K}
\\
\mathrm{V}(\mathbb{Q}_{\lambda_0,r_0},\mathbb{Q}_{\lambda,r}) & \leq2
\lambda_0(1+\lambda_0)\mathrm{V}(\mathbb{P}_{r_0},
\mathbb {P}_{r})+ 4\lambda_0\mathrm{K}(
\mathbb{P}_{r_0},\mathbb{P}_{r})\nonumber
\\
& \quad+2\mathrm{V}(\lambda_0,\lambda)+4\mathrm{K}(
\lambda_0,\lambda )+2\mathrm{K}(\lambda_0,
\lambda)^2,
\label{eq:V}
\\
h(\mathbb{Q}_{\lambda_0,r_0},\mathbb{Q}_{\lambda,r}) &\leq \sqrt {
\lambda_0}\,h(\mathbb{P}_{r_0},\mathbb{P}_{r})+
h(\lambda_0,\lambda ).\label{eq:h}
\end{align}
Moreover, there exists a constant $\overline{C}\in(0,\infty)$,
depending on $\underline{\lambda}$ and $\overline{\lambda}$ only, such
that for all $\lambda_0,\lambda\in[\underline{\lambda},\overline
{\lambda}]$,
\begin{align}
\mathrm{K}(\mathbb{Q}_{\lambda_0,r_0},\mathbb{Q}_{\lambda,r}) & \leq
\overline{C} \bigl(\mathrm{K}(\mathbb{P}_{r_0},\mathbb{P}_{r})+|
\lambda _0-\lambda|^2\bigr),\label{eq:1K}
\\
\mathrm{V}(\mathbb{Q}_{\lambda_0,r_0},\mathbb{Q}_{\lambda,r}) & \leq
\overline{C} \bigl(\mathrm{V}(\mathbb{P}_{r_0},\mathbb{P}_{r})+
\mathrm {K}(\mathbb{P}_{r_0},\mathbb{P}_{r})+|
\lambda_0-\lambda|^2\bigr),\label
{eq:1V}
\\
h(\mathbb{Q}_{\lambda_0,r_0},\mathbb{Q}_{\lambda,r})&\leq\overline{C} \bigl(|
\lambda_0-\lambda| + h(\mathbb{P}_{r_0},
\mathbb{P}_{r})\bigr).\label{eq:1h} 
\end{align}
%
\end{lemma}
The proof of the lemma is given in Section \ref{pr.lem.1}. We proceed with
the proof of Theorem \ref{mainthm}.

Let $\varepsilon_n=n^{-\gamma}(\log n)^{\ell}$ for $\gamma$ and $\ell
>\ell_0$ as in the statement of Theorem \ref{mainthm}. Set $\overline
{\varepsilon}_n=2\overline{C}\varepsilon_n$, where $\overline{C}$ is
the constant from Lemma \ref{lem:ineq}. We define the sieves of
densities $\mathcal{F}_n$ as in Theorem 5 in \cite{shen13}:
\begin{align*}
\mathcal{F}_n= \Biggl\{ r_{F,\varSigma} \textrm{ with } F=\sum
_{i=1}^{\infty
} \pi_i
\delta_{z_i}:{} & z_i\in[-\alpha_n,
\alpha_n]^d,\forall i\leq I_n; \sum
_{i>I_n}\pi_i <\varepsilon_n;
\\
&\sigma_{0,n}^2 \leq\operatorname{eig}_j(
\varSigma)<\sigma _{0,n}^2\bigl(1+\varepsilon_n^2/d
\bigr)^{J_n} \Biggr\},
\end{align*}
where
\begin{equation*}
I_n=\bigl\lfloor n\varepsilon_n^2 /\log n
\bigr\rfloor, \qquad J_n=\alpha _n^{a_1}=
\sigma_{0,n}^{-2a_2}=n,
\end{equation*}
and $a_1$ and $a_2$ are as in Assumption \ref{ass:prior}.
We also put
%
\begin{equation}
\label{Qn} \mathcal{Q}_n=\bigl\{ \mathbb{Q}_{\lambda,r}:r\in
\mathcal{F}_n,\lambda\in [\underline{\lambda},\overline{\lambda}]
\bigr\}.
\end{equation}

In \cite{shen13}, sieves of the type $\mathcal{F}_n$ are used to verify
conditions of Theorem \ref{thm2.1ghosal01} and to determine posterior
contraction rates in the standard density estimation context. We
further will show that these sieves also work in the case of
decompounding by verifying the conditions of Theorem \ref
{thm2.1ghosal01} for the sieves $\mathcal{Q}_n$ defined in \eqref{Qn}.
\endgroup

\subsection{Verification of \eqref{c1}}

Introduce the notation
\begin{equation*}
\overline{h}_1(\lambda_1,\lambda_2)=
\overline{C}|\lambda_1-\lambda_2|, \qquad
\overline{h}_2(r_1,r_2)=\overline{C}h(
\mathbb{P}_{r_1},\mathbb{P}_{r_2}).
\end{equation*}
Let $\{\lambda_i\}$ be the centers of the balls from a minimal covering
of $[\underline{\lambda}, \overline{\lambda}]$ with $\overline
{h}_1$-intervals of size $\overline{C}\varepsilon_n$. Let $\{r_j\}$ be
centers of the balls from a minimal covering of $\scr{F}_{n}$ with
$\overline{h}_2$-balls of size $\overline{C}\varepsilon_n$. By Lemma
\ref{lem:ineq}, for any $\bb{Q}_{\lambda, r} \in\scr{Q }_{n}$,
\[
h(\bb{Q}_{\lambda,r} , \bb{Q}_{\lambda_i, r_j}) \le\overline {h}_1(
\lambda, \lambda_i) + \overline{h}_2(r, r_j)
\le\overline {\varepsilon}_n
\]
by appropriate choices of $i$ and $j$. Hence,
\begin{equation*}
N (\overline{\varepsilon}_n,\mathcal{Q}_{n},{h}) \leq N
\bigl(\overline {C}\varepsilon_n,[\underline{\lambda},\overline{
\lambda}],\overline {h}_1\bigr) \times N ( \overline{C}
\varepsilon_n , {\mathcal{F}}_n , \overline{h}_2
),
\end{equation*}
and so
\begin{equation*}
\log N(\overline{\varepsilon}_n,\mathcal{Q}_n,h)\leq\log
N\bigl(\overline {C}\varepsilon_n,[\underline{\lambda},\overline{
\lambda}],\overline {h}_1\bigr)+\log N(\overline{C}
\varepsilon_n,\mathcal{F}_n,\overline{h}_2).
\end{equation*}
By Proposition 2 and Theorem 5 in \cite{shen13}, there exists a
constant $c_1>0$ such that for all $n$ large enough,
\begin{equation*}
\log N(\overline{C}\varepsilon_n,\mathcal{F}_n,
\overline{h}_2) = \log N(\varepsilon_n,
\mathcal{F}_n,h)\leq c_1 n\varepsilon_n^2=
\frac
{c_1}{4\overline{C}^2}n\overline{\varepsilon}_n^2.
\end{equation*}
On the other hand,
\begin{align*}
\log N\bigl(\overline{C}\varepsilon_n,[\underline{\lambda},
\overline{\lambda }],\overline{h}_1\bigr) & = \log N\bigl(
\varepsilon_n,[\underline{\lambda },\overline{\lambda}],|\cdot|\bigr),
\\
&\lesssim\log \biggl( \frac{1}{\varepsilon_n} \biggr)
\\
&\lesssim\log \biggl( \frac{1}{\overline{\varepsilon}_n} \biggr).
\end{align*}
With our choice of $\overline{\varepsilon}_n$, for all $n$ large
enough, we have
\begin{equation*}
\frac{c_1}{4\overline{C}^2}n\overline{\varepsilon}_n^2 \geq\log
\biggl( \frac{1}{\overline{\varepsilon}_n} \biggr),\vadjust{\goodbreak}
\end{equation*}
so that for all $n$ large enough,
\begin{equation*}
\log N(\overline{\varepsilon}_n,\mathcal{Q}_n,h)\leq
\frac
{c_1}{2\overline{C}^2} n\overline{\varepsilon}_n^2.
\end{equation*}
We can simply rename the constant $c_1/(2\overline{C}^2)$ in this
formula into $c_1$, and thus \eqref{c1} is satisfied with that constant.

\subsection{Verification of \eqref{c2} and \eqref{c3}}

We first focus on \eqref{c3}.
Introduce
\begin{equation*}
\widetilde{B}(\eps, \bb{Q}_{\lambda_0, r_0})= \bigl\{ (\lambda,r): \mathrm {K}(
\mathbb{P}_{r_0},\mathbb{P}_{r}) \leq\varepsilon^{2},
\mathrm {V}(\mathbb{P}_{r_0},\mathbb{P}_{r}) \leq
\varepsilon^{2}, |\lambda _0-\lambda|\leq\varepsilon\bigr
\}.
\end{equation*}
Suppose that $(\lambda, r) \in\tilde{B}(\eps, \bb{Q}_{\lambda_0,
r_0})$. From \eqref{eq:1K} we obtain
\[
\mathrm{K}(\bb{Q}_{\lambda_0, r_0}, \bb{Q}_{\lambda, r}) \le\overline {C}
\mathrm{K}(\bb{P}_{r_0}, \bb{P}_r) + \overline{C} |\lambda-
\lambda _0|^2 \le2\overline{C} \eps^2.
\]
Furthermore, using \eqref{eq:1V}, we have
\begin{align*}
\mathrm{V}(\bb{Q}_{\lambda_0, r_0}, \bb{Q}_{\lambda, r}) &\le\overline{C}
\mathrm{V}(\bb{P}_{r_0}, \bb{P}_r) + \overline{C}
\mathrm{K}(\bb{P}_{r_0}, \bb{P}_r) +\overline{C} |\lambda-
\lambda_0|^2 %
\le3\overline{C}
\eps^2.
\end{align*}
Combination of these inequalities with the definition of the set
$B(\varepsilon,\mathbb{Q}_{\lambda_0,r_0})$ in \eqref{set_b} yields
\begin{equation*}
\widetilde{B}(\eps, \mathbb{Q}_{\lambda_0,r_0} ) \subset B( \sqrt {3\overline{C}}
\eps,\mathbb{Q}_{\lambda_0,r_0} ).
\end{equation*}
Consequently,
\begin{align}
\varPi \bigl( B( \sqrt{3\overline{C}} \eps,\mathbb{Q}_{\lambda_0,r_0}) \bigr)
&{}\geq\varPi \bigl( \tilde{B}(\eps, \mathbb{Q}_{\lambda_0,r_0} ) \bigr)\nonumber\\
&{}=\varPi_1 ( |\lambda_0-\lambda| \leq\eps )\nonumber\\
&\quad{}\times\varPi_2 \bigl( r_{f,\varSigma}: \mathrm{K}(
\pp_{r_0},\pp _{r_{F,\varSigma}})\leq\eps^2,
\, \mathrm{V}(\pp_{r_0},\pp_{r_{F,\varSigma}}) \leq
\eps^2 \bigr).\label{eq:pipi}
\end{align}
By Assumption \ref{ass:prior}(i),
\begin{equation*}
\varPi_1 ( |\lambda_0-\lambda| \leq\eps ) \geq
\underline{\pi }_1 \eps.
\end{equation*}
%
Furthermore,
Theorem 4 in \cite{shen13} yields that for some $A,C>0$ and all
sufficiently large $n$,
\begin{gather*}
\varPi_2 \bigl( r_{F,\varSigma}: \mathrm{K}(\pp_{r_0},
\pp_{r_{F,\varSigma}})\leq {A} n^{-2\gamma} (\log n)^{2\ell_0}, \mathrm{V}(
\pp_{r_0},\pp _{r_{F,\varSigma}}) \leq{A} n^{-2\gamma} (\log
n)^{2\ell_0} \bigr)
\\
\geq\exp \bigl( -C n \bigl\{ n^{-\gamma} (\log n)^{\ell_0} \bigr
\}^2 \bigr).
\end{gather*}
We substitute $\eps$ with $\sqrt{A}n^{-\gamma}(\log n)^{\ell_0}$ and
write $\widetilde{\eps}_n=\sqrt{3A\overline{C}}n^{-\gamma}(\log n)^{\ell
_0}$ to arrive at
\[
\varPi \bigl( B( \widetilde{\eps}_n,\mathbb{Q}_{\lambda_0,r_0})
\bigr)\geq \underline{\pi}_1 \sqrt{A}n^{-\gamma}(\log
n)^{\ell_0}\times\exp \biggl( -\frac{C}{3A\overline{C} }n \widetilde{
\eps}_n^2 \biggr).
\]
Now, since $\gamma<\frac{1}{2}$,
for all $n$ large enough, we have
\begin{equation*}
\underline{\pi}_1 \sqrt{A} n^{-\gamma} (\log n)^{\ell_0}
\geq\exp \bigl( - n^{1-2\gamma} (\log n)^{2\ell_0} \bigr).
\end{equation*}
Consequently, for all $n$ large enough,
%
\begin{equation}
\label{c3_ineq} \varPi ( B( \widetilde{\varepsilon}_n ,
\mathbb{Q}_{\lambda_0,f_0} ) \geq\exp \biggl( - \biggl(\frac{C+1}{3A\overline{C}} \biggr) n
\widetilde{\eps}_n^2 \biggr).
\end{equation}
Choosing $c_2=\frac{C+1}{3A\overline{C}}$, we have verified \eqref{c3}
(with $c_4=1$).

For the verification of \eqref{c2}, we use the constants $c_2$ and
$\widetilde\eps_n$ as above. Note first that
\begin{equation*}
\varPi(\mathcal{Q}\setminus\mathcal{Q}_n)=\varPi_2\bigl(
\mathcal{F}_n^c\bigr).
\end{equation*}
By Theorem 5 in \cite{shen13} (see\ also p.~627 there), for some
$c_3>0$ and any constant $c>0$, we have
\begin{equation*}
\varPi_2\bigl(\mathcal{F}_n^c\bigr)\leq
c_3 \exp \bigl(-(c+4)n \bigl\{ n^{-\gamma} (\log
n)^{\ell_0} \bigr\}^2 \bigr),
\end{equation*}
provided that $n$ is large enough. Thus,
%
\begin{equation*}
\varPi(\mathcal{Q}\setminus\mathcal{Q}_n) \leq c_3 \exp
\biggl(-\frac
{c+4}{3A\overline{C}}\, n \widetilde{\varepsilon}_n^2
\biggr).
\end{equation*}
Without loss of generality, we can take the positive constant $c$
greater than $3A\overline{C}(c_2+4)-4$. This gives
\begin{equation*}
\varPi(\mathcal{Q}\setminus\mathcal{Q}_n)\leq c_3\exp
\bigl(-(c_2+4) n \widetilde{\varepsilon}_n^2
\bigr),
\end{equation*}
which is indeed \eqref{c2}.

We have thus verified conditions \eqref{c1}--\eqref{c3}, and the
statement of Theorem \ref{mainthm} follows by Theorem~\ref
{thm2.1ghosal01} since $\bar\eps_n\geq\widetilde\eps_n$ (eventually).

\section{Proof of Lemma \ref{lem:ineq}}
\label{pr.lem.1}

We start with a lemma from \cite{csiszar1963}, which will be used three
times in the proof of Lemma~\ref{lem:ineq}.
%
Consider a probability space $(\varOmega,\mathfrak{F},\pp)$. Let $\pp_0$
be a probability measure on $(\varOmega,\mathfrak{F})$ and assume that $\pp
_0\ll\pp$ with Radon--Nikodym derivative $\zeta=\frac{\dd\pp_0}{\dd
\pp}$. Furthermore, let $\mathfrak{G}$ be a sub-$\sigma$-algebra of
$\mathfrak{F}$. The restrictions of $\pp$ and $\pp_0$ to $\mathfrak{G}$
are denoted $\pp'$ and $\pp'_0$, respectively. Then $\pp'_0\ll\pp'$
and $\frac{\dd\pp'_0}{\dd\pp'}=\ee_{\pp}[\zeta|\mathfrak{G}]=:\zeta'$.

\begin{lemma}\label{lem:convexnew}
Let $g:[0,\infty)\to\RR$ be a convex function. Then
\[
\ee_{\pp'}g\bigl(\zeta'\bigr)\leq\ee_\pp\, g(
\zeta).
\]
\end{lemma}
%
The proof of the lemma consists in an application of Jensen's
inequality for conditional expectations. This lemma is typically used
as follows. The measures $\pp$ and $\pp_0$ are possible distributions
of some random element $X$. If $X'=T(X)$ is some measurable
transformation of $X$, then we consider $\pp'$ and $\pp'_0$ as the
corresponding distributions of $X'$. Here $T$ may be a projection. In
the present context, we take $X=(X_t, t\in[0,1])$ and $X'=X_1$, and so
$\pp$ in the lemma should be taken as $\RR=\RR_{\lambda,r}$ and $\pp'$
as $\mathbb{Q}=\mathbb{Q}_{\lambda,r}$.

In the proof of Lemma \ref{lem:ineq}, for economy of notation, a
constant $c(\underline{\lambda},\overline{\lambda})$ depending on
$\underline{\lambda}$ and $\overline{\lambda}$ may differ from line to
line. We also abbreviate $\mathbb{Q}_{\lambda_0,r_0}$ and $\mathbb
{Q}_{\lambda,r}$ to $\mathbb{Q}_{0}$ and $\mathbb{Q}$, respectively.
The same convention will be used for $\bb{R}_{\lambda_0, r_0}$, $\bb
{R}_{\lambda, r}$, $\bb{P}_{r_0}$, and $\bb{P}_{r}$.

\begin{proof}[Proof of inequalities~(\ref{eq:K}) and~(\ref{eq:1K})]
Application of Lemma~\ref{lem:convexnew} with $g(x)=\break (x\log x)
1_{\{ x\geq0\}}$ gives $ \mathrm{K}(\mathbb{Q}_{0},\mathbb{Q})
\leq\mathrm{K}(\mathbb {R}_{0},\mathbb{R})$. Using \eqref{eq:contlik}
and the expression for the mean of a stochastic integral with respect
to a Poisson point process (see, e.g.,\ property 6 on p.~68
in~\cite{skorohod64}), we obtain that
\begin{align*}
\mathrm{K}(\mathbb{R}_{0},\mathbb{R}) & = \int\log \biggl(
\frac{\mathrm
{d}\mathbb{R}_0}{\mathrm{d}\mathbb{R}} \biggr)\mathrm{d}\mathbb{R}_0
\\
&=\lambda_0\int\log \biggl( \frac{\lambda_0 r_0}{\lambda r} \biggr)
r_0 - (\lambda_0-\lambda)
\\
&= \lambda_0 \mathrm{K}(\mathbb{P}_{0},\mathbb{P})+
\biggl( \lambda_0 \log \biggl( \frac{\lambda_0}{\lambda} \biggr) - [
\lambda_0-\lambda] \biggr)
\\
&= \lambda_0 \mathrm{K}(\mathbb{P}_{0},\mathbb{P})+
\mathrm{K}(\lambda _0,\lambda).
\end{align*}
Now
\begin{align*}
\lambda_0 \log \biggl( \frac{\lambda_0}{\lambda} \biggr) - (\lambda
_0-\lambda) & = \lambda_0 \biggl\llvert \log \biggl(
\frac{\lambda}{\lambda_0} \biggr) - \biggl( \frac
{\lambda}{\lambda_0} -1 \biggr) \biggr\rrvert
\\
& \leq c(\underline{\lambda},\overline{\lambda})|\lambda_0-
\lambda|^2,
\end{align*}
where $c(\underline{\lambda},\overline{\lambda})$ is some constant
depending on $\underline{\lambda}$ and $\overline{\lambda}$. The result follows.
\end{proof}

%
\begin{proof}[Proof of inequalities~(\ref{eq:V}) and~(\ref{eq:1V})]
We have
\begin{align*}
\mathrm{V}(\mathbb{Q}_0,\mathbb{Q}) & = \ee_{\mathbb{Q}_0} \biggl[
\log^2 \biggl( \frac{\mathrm{d}\mathbb
{Q}_0}{\mathrm{d}\mathbb{Q}} \biggr) 1_{ \{ \frac{\mathrm{d}\mathbb
{Q}_0}{\mathrm{d}\mathbb{Q}} \geq1  \}} \biggr]
+ \ee_{\mathbb{Q}_0} \biggl[ \log^2 \biggl( \frac{\mathrm{d}\mathbb
{Q}_0}{\mathrm{d}\mathbb{Q}}
\biggr) 1_{ \{ \frac{\mathrm{d}\bb
{Q}_0}{\dd\bb{Q}} < 1  \}} \biggr]
\\
&=\mathrm{I}+\mathrm{II}.
\end{align*}
Application of Lemma~\ref{lem:convexnew} with $g(x)=(x\log^2(x)) 1_{\{
x\geq1\}} $ (which is a convex function) gives
%
\begin{equation}
\label{ineqI} \mathrm{I} \leq \ex_{\bb{R}_0} \biggl[ \log^2
\biggl( \frac{ \mathrm
{d}\mathbb{R}_{0} }{ \mathrm{d}\mathbb{R} } \biggr) 1_{  [ \frac{
\mathrm{d}\mathbb{R}_{0} }{ \mathrm{d}\mathbb{R} } \geq1  ] } \biggr] \le\mathrm{V}(
\bb{R}_0, \bb{R}).
\end{equation}
As far as $\mathrm{II}$ is concerned, for $x\ge0$, we have the inequalities
\[
\frac{x^2}{2} \le e^x-1-x \le2 \bigl(e^{x/2}-1
\bigr)^2.
\]
The first inequality is trivial, and the second is a particular case of
inequality~(8.5) in \xch{\cite{ghosal00}}{Ghosal, Ghosh, and Van der Vaart (2000)} and is
equally elementary.
The two inequalities together yield
\[
e^{-x}x^2 \le4 \bigl(e^{-x/2}-1
\bigr)^2.
\]
Applying this inequality with $x=-\log\frac{\dd\bb{Q}_0}{\dd\bb{Q}}$
(which is positive on the event \mbox{$\{\frac{\dd\bb{Q}_0}{\dd\bb{Q}}
{<} 1\} $}) and taking the expectation with respect to $\bb{Q}$ give
\begin{align*}
\mathrm{II} &= \ee_{\bb{Q}} \biggl[ \frac{\dd\bb
{Q}_0}{\dd\bb{Q}}\log^2
\frac{\dd\bb{Q}_0}{\dd\bb{Q}} 1_{\{ \frac{\dd
\bb{Q}_0}{\dd\bb{Q}}< 1\}} \biggr]
\\
& \le4 \int \biggl( \sqrt{\frac{\dd\bb{Q}_0}{\dd\bb{Q}}} -1 \biggr)^2 \dd\bb{Q}
\\
&= 4 h^2(\bb{Q}_0, \bb{Q}) \le4\mathrm{K}(
\bb{Q}_0, \bb{Q}). 
\end{align*}
For the final inequality, see \cite{pollard2002}, p.~62, formula~(12).

Combining the estimates on $\mathrm{I}$ and $\mathrm{II,}$ we obtain that
%
\begin{equation}
\label{VQP} \mathrm{V}(\mathbb{Q}_0,\mathbb{Q}) \leq\mathrm{V}(
\mathbb {R}_{0},\mathbb{R}) + 4 \mathrm{K}(\bb{Q}_0,
\bb{Q}).
\end{equation}
After some long and tedious calculations employing \eqref{eq:contlik}
and the expressions for the mean and variance of a stochastic integral
with respect to a Poisson point process (see, e.g.,\ property 6 on
p.~68 in~\cite{skorohod64} and Lemma 1.1 in \cite{kutoyants98}), we get that
\begin{align*}
\mathrm{V}(\mathbb{R}_0,\mathbb{R}) & = \lambda_0 \int
\biggl\{ \log \biggl( \frac{\lambda_0}{\lambda} \biggr) + \log \biggl(
\frac{r_0}{r} \biggr) \biggr\}^2 f_0
\\
&\quad+ \lambda_0^2 \biggl\{ \int\log \biggl(
\frac{r_0}{r} \biggr) r_0 + \log \biggl( \frac{\lambda_0}{\lambda}
\biggr) - \biggl( 1 - \frac
{\lambda}{\lambda_0} \biggr) \biggr\}^2
\\
&=\mathrm{III}+\mathrm{IV}.
\end{align*}
By the $c_2$-inequality $(a+b)^2\leq2a^2+2b^2$ we have
%
\begin{align}
\mathrm{III}
&\leq2 \lambda_0 \log^2 \biggl( \frac{\lambda_0}{\lambda} \biggr) + 2\lambda_0 \int\log^2 \biggl(\frac{r_0}{r} \biggr)r_0\nonumber\\
& = 2\mathrm{V}(\lambda_0,\lambda) + 2\lambda_0\mathrm{V}(\mathbb {P}_{0},\mathbb{P}),\label{ineqIII}
\end{align}
from which we deduce
%
\begin{equation}
\label{upIII} \mathrm{III}\leq c( \underline{\lambda},\overline{\lambda} ) |
\lambda _0-\lambda|^2 + 2\overline{\lambda} \mathrm{V}(
\mathbb{P}_{0},\mathbb{P})
\end{equation}
for some constant $c( \underline{\lambda},\overline{\lambda} )$
depending on $\underline{\lambda}$ and $\overline{\lambda}$ only. As
far as $\mathrm{IV}$ is concerned, the $c_2$-inequality and the
Cauchy--Schwarz inequality give that
%
\begin{align}
\mathrm{IV} & \leq2
\lambda_0^2 \biggl( \int\log \biggl( \frac{r_0}{r}
\biggr)r_0 \biggr)^2+2\lambda_0^2
\biggl( \log \biggl(\frac{\lambda
_0}{\lambda} \biggr) - \biggl[ 1 - \frac{\lambda}{\lambda_0}
\biggr] \biggr)^2\nonumber
\\
& \leq2 \lambda_0^2 \mathrm{V}(\mathbb{P}_{0},
\mathbb{P})+2\mathrm {K}(\lambda_0,\lambda)^2,\label{ineqIV}
\end{align}
from which we find the upper bound
%
\begin{equation}
\label{upIV} \mathrm{IV} \leq2 \overline{\lambda}^2 \mathrm{V}(
\mathbb {P}_{0},\mathbb{P})+c(\underline{\lambda},\overline{\lambda}) |
\lambda _0-\lambda|^2
\end{equation}
for some constant $c(\underline{\lambda},\overline{\lambda})$ depending
on $\underline{\lambda}$ and $\overline{\lambda}$.
Combining estimates \eqref{ineqIII} and \eqref{ineqIV} on $\mathrm
{III}$ and $\mathrm{IV}$ with inequalities \eqref{VQP} and \eqref{eq:K}
yields \eqref{eq:V}. Similarly, the upper bounds \eqref{upIII} and
\eqref{upIV}, combined with \eqref{VQP} and \eqref{eq:K}, yield \eqref{eq:1V}.
%
\end{proof}

\begin{proof}[Proof of inequalities (\ref{eq:h}) and (\ref{eq:1h})]
First, note that for $g(x) = (\sqrt{x}-1)^2 1_{[x\geq0]}$,
\[
h^2(\bb{Q}_0, \bb{Q}) = \ee_{\bb{Q}} \biggl[
\biggl(\sqrt{\frac{\dd\bb
{Q}_0}{\dd\bb{Q}}}-1 \biggr)^2 \biggr] =
\ee_{\bb{Q}} \biggl[ g \biggl( \frac{\dd\bb{Q}_0}{\dd\bb{Q}} \biggr) \biggr].
\]
Since $g$ is convex, an application of Lemma~\ref{lem:convexnew} yields
$h(\mathbb{Q}_{0},\mathbb{Q}) \leq h(\mathbb{R}_{0},\mathbb{R})$.
Using \eqref{eq:contlik} and invoking Lemma 1.5 in~\cite{kutoyants98},
in particular, using formula (1.30) in its statement, we get that
%
\begin{align*}
h(\mathbb{R}_{0},\mathbb{R}) & \leq\| \sqrt{\lambda_0
r_0} - \sqrt {\lambda r}\|
\\
& \leq\| \sqrt{\lambda_0 r_0} - \sqrt{
\lambda_0 r}\|+\| \sqrt{\lambda _0 r} - \sqrt{\lambda r}
\|
\\
& \leq\sqrt{\lambda_0} \| \sqrt{ r_0} - \sqrt{r}\|+|
\sqrt{\lambda_0} - \sqrt{\lambda}|
\\
& = \sqrt{\lambda_0} h(\pp_{0},\pp) + h(
\lambda_0,\lambda),
\end{align*}
where $\|\cdot\|$ denotes the $L^2$-norm. This proves \eqref{eq:h}.
Furthermore, from this we obtain the obvious upper bound
\[
h(\mathbb{R}_{0},\mathbb{R}) \leq \sqrt{\overline{\lambda}}\, h(\pp
_{0},\pp) + \frac{1}{2\sqrt{\underline{\lambda}}}|\lambda_0-\lambda|,
\]
which yields \eqref{eq:1h}.
\end{proof}

\section*{Acknowledgments}
The authors would like to thank the referee for his/her remarks. The
research leading to these results has received funding from the
European Research Council under ERC Grant Agreement 320637.

%

\end{document}